# Recurrent extensions of self-similar Markov processes and Cramér's condition II

VíCTOR RIVERO

[1]*CIMAT A.C. Calle Jalisco s/n, Col. Valenciana, C.P. 36240 Guanajuato, Guanajuato, Mexico.
E-mail: rivero@cimat.mx*

We prove that a positive self-similar Markov process $(X, \mathbb{P})$ that hits 0 in a finite time admits a self-similar recurrent extension that leaves 0 continuously if and only if the underlying Lévy process satisfies Cramér's condition.

*Keywords:* excursion theory; exponential functionals of Lévy processes; Lamperti's transformation; Lévy processes; self-similar Markov processes

## 1. Introduction and main result

Let $\mathbb{P} = (\mathbb{P}_x, x \geq 0)$ be a family of probability measures on Skohorod's space $\mathbb{D}^+$, the space of càdlàg paths defined on $[0, \infty[$ with values in $\mathbb{R}^+$. The space $\mathbb{D}^+$ is endowed with the Skohorod topology and its Borel $\sigma$-field. We will denote by $X$ the canonical process of the coordinates and $(\mathcal{G}_t, t \geq 0)$ will be the natural filtration generated by $X$. Assume that under $\mathbb{P}$ the canonical process $X$ is a positive self-similar Markov process (pssMp). That is, $(X, \mathbb{P})$ is a $[0, \infty[$-valued strong Markov process with the following scaling property: there exists an $\alpha > 0$ such that for every $c > 0$,

$$(\{cX_{tc^{-1/\alpha}}, t \geq 0\}, \mathbb{P}_x) \stackrel{\text{Law}}{=} (\{X_t, t \geq 0\}, \mathbb{P}_{cx}) \qquad \forall x \geq 0.$$

We will further assume that $(X, \mathbb{P})$ is a pssMp that hits 0 in a $\mathbb{P}$-a.s. finite time $T_0 = \inf\{t > 0 : X_t = 0\}$ and dies. So, $\mathbb{P}_0$ is the law of the degenerate path equal to 0. According to Lamperti's transformation [14], the family of laws $\mathbb{P}$ can be obtained as the image law of the exponential of an $\mathbb{R} \cup \{-\infty\}$-valued Lévy process $\xi$ with law $\mathbf{P}$, time-changed by the inverse of the additive functional,

$$t \to \int_0^t \exp\{\xi_s/\alpha\} \, \mathrm{d}s, \qquad t \geq 0. \tag{1}$$

As usual, any function $f : \mathbb{R} \to \mathbb{R}$ is extended to $\mathbb{R} \cup \{-\infty\}$ by taking $f(-\infty) = 0$. Thus, the state $\{-\infty\}$ will be taken as a cemetery state for $\xi$ and we denote by $\zeta$ its lifetime, namely $\zeta := \inf\{t > 0 : \xi_t = -\infty\}$, and by $\{\mathcal{F}_t, t \geq 0\}$, the filtration of $\xi$. A consequence of Lamperti's transformation is that the law of $T_0$ under $\mathbb{P}_x$ is equal to that of $x^{1/\alpha} I$ under $\mathbf{P}$, where $I$ denotes the exponential functional associated to $\xi$,

$$I := \int_0^\zeta \exp\{\xi_s/\alpha\}\, ds.$$

Lamperti proved the following characterization for pssMp that hit 0 in a finite time: either $(X, \mathbb{P})$ hits 0 by a jump and in a finite time

$$\mathbb{P}_x(T_0 < \infty, X_{T_0-} > 0, X_{T_0+t} = 0,\ \forall\, t \geq 0) = 1 \qquad \forall x > 0,$$

which happens if and only if $\mathbf{P}(\zeta < \infty) = 1$, or $(X, \mathbb{P})$ hits 0 continuously and in a finite time

$$\mathbb{P}_x(T_0 < \infty, X_{T_0-} = 0, X_{T_0+t} = 0,\ \forall\, t \geq 0) = 1 \qquad \forall x > 0,$$

and this is equivalent to $\mathbf{P}(\zeta = \infty, \lim_{t \to \infty} \xi_t = -\infty) = 1$. Reciprocally, the image law of the exponential of any $\mathbb{R} \cup \{-\infty\}$-valued Lévy process time-changed by the inverse of the functional defined in (1) is the law of a pssMp that dies at its first hitting time of 0. For more details, see [14] or [16].

The main purpose of this note is to continue our study, initiated in [16], on the existence and characterization of positive-valued self-similar Markov processes $\widetilde{X}$ that behave like $(X, \mathbb{P})$ before their first hitting time of 0 and for which the state 0 is a regular and recurrent state. Such a process $\widetilde{X}$ will be called a *recurrent extension* of $(X, \mathbb{P})$. We refer to [16, 17] and the references therein for an introduction to this problem and for background on excursion theory for positive self-similar Markov processes.

We say that a $\sigma$-finite measure n on $(\mathbb{D}^+, \mathcal{G}_\infty)$ having infinite mass is an *excursion measure* compatible with $(X, \mathbb{P})$ if the following are satisfied:

(i) n is carried by

$$\{\omega \in \mathbb{D}^+ | T_0(\omega) > 0 \text{ and } X_t(\omega) = 0, \forall t \geq T_0\};$$

(ii) for every bounded $\mathcal{G}_\infty$-measurable $H$ and each $t > 0$ and $\Lambda \in \mathcal{G}_t$,

$$\mathrm{n}(H \circ s_t, \Lambda \cap \{t < T_0\}) = \mathrm{n}(\mathbb{E}_{X_t}(H), \Lambda \cap \{t < T_0\}),$$

where $s_t$ denotes the shift operator;

(iii) $\mathrm{n}(1 - \mathrm{e}^{-T_0}) < \infty$.

Moreover, we will say that n is *self-similar* if it has the following scaling property: there exists a $0 < \gamma < 1$, s.t. for all $a > 0$, the measure $H_a \mathrm{n}$, which is the image of n under the mapping $H_a : \mathbb{D}^+ \to \mathbb{D}^+$, defined by

$$H_a(\omega)(t) = a\omega(a^{-1/\alpha}t), \qquad t \geq 0,$$

is such that

$$H_a \mathrm{n} = a^{\gamma/\alpha} \mathrm{n}.$$

The parameter $\gamma$ will be called the *index of self-similarity* of n. See Section 2 in [16] for equivalent definitions of self-similar excursion measure.

The entrance law associated with n is the family of finite measures $(\mathrm{n}_t, t > 0)$, defined by

$$\mathrm{n}(X_t \in \mathrm{d}y, t < T_0) = \mathrm{n}_t(\mathrm{d}y), \qquad t > 0.$$

It is known that there exists a one to one correspondence between recurrent extensions of $(X, \mathbb{P})$ and self-similar excursion measures compatible with $(X, \mathbb{P})$; see, for example, [16, 17]. So, determining the existence of the recurrent extensions of $(X, \mathbb{P})$ is equivalent to doing so for self-similar excursion measures. We recall that the index of self-similarity of a self-similar excursion measure coincides with that of the stable subordinator which is the inverse of the local time at 0 of the associated recurrent extension of $(X, \mathbb{P})$.

We say that a positive self-similar Markov process for which 0 is a regular and recurrent state *leaves* 0 *continuously* (resp., *by a jump*) whenever its excursion measure n is carried by the paths that leave 0 continuously (resp., that leave 0 by a jump)

$$\mathrm{n}(X_{0+} > 0) = 0 \qquad [\text{resp., } \mathrm{n}(X_{0+} = 0) = 0].$$

Vuolle-Apiala [17] proved, under some hypotheses, that any positive self-similar Markov process for which 0 is a regular and recurrent state either leaves 0 continuously or by jumps. In fact, his result still holds true in the general setting, as is proved in the following lemma.

**Lemma 1.** *Let* n *be a self-similar excursion measure compatible with* $(X, \mathbb{P})$ *and with index of self-similarity* $\gamma \in \,]0,1[$. *Then,*

$$\textit{either } \mathrm{n}(X_{0+} > 0) = 0 \quad \textit{or} \quad \mathrm{n}(X_{0+} = 0) = 0.$$

**Proof.** Assume that the claim of the lemma does not hold. Let $\mathrm{n}^c = c^{(c)} \mathrm{n}|_{\{X_{0+}=0\}}$ and $\mathrm{n}^j = c^{(j)} \mathrm{n}|_{\{X_{0+}>0\}}$ be the restrictions of n to the set of paths $\{X_{0+} = 0\}$ and $\{X_{0+} > 0\}$, respectively, and $c^{(c)}$ and $c^{(j)}$ be normalizing constants such that

$$\mathrm{n}^c(1 - \mathrm{e}^{-T_0}) = 1 = \mathrm{n}^j(1 - \mathrm{e}^{-T_0}).$$

The measures $\mathrm{n}^c$ and $\mathrm{n}^j$ are self-similar excursion measures compatible with $(X, \mathbb{P})$ and with the same self-similarity index $\gamma$. According to Lemma 3 of [16], the potential measure of $\mathrm{n}^c$ and that of $\mathrm{n}^j$ are given by the same purely excessive measure

$$\mathrm{n}^c\left(\int_0^{T_0} 1_{\{X_t \in \mathrm{d}y\}} \, \mathrm{d}t\right) = C_{\alpha,\gamma} y^{(1-\alpha-\gamma)/\alpha} \, \mathrm{d}y = \mathrm{n}^j\left(\int_0^{T_0} 1_{\{X_t \in \mathrm{d}y\}} \, \mathrm{d}t\right), \qquad y > 0, \quad (2)$$

where $C_{\alpha,\gamma} \in \,]0, \infty[$ is a constant. So, by Theorem 5.25 of [11] on the uniqueness of purely excessive measures, the entrance laws associated with $\mathrm{n}^c$ and $\mathrm{n}^j$ are equal. Hence, by Theorem 4.7 of [6], the measures $n^c$ and $n^j$ are equal. This leads to a contradiction of the fact that the supports of the measures $n^c$ and $n^j$ are disjoint. □

If $n^\beta$ is a self-similar excursion measure with index $\gamma = \beta\alpha \in\ ]0,1[$ and is carried by the paths that leave 0 by a jump, then the self-similarity implies that $n^\beta$ has the form $n^\beta = c_{\alpha,\beta}\mathbb{P}_{\eta_\beta}$, where $0 < c_{\alpha,\beta} < \infty$ is a normalizing constant and the starting measure or jumping-in measure $\eta_\beta$ is given by

$$\eta_\beta(\mathrm{d}x) = \beta x^{-1-\beta}\,\mathrm{d}x, \qquad x > 0.$$

The choice of the constant $c_{\alpha,\beta}$ depends on the normalization of the local time at 0 of the recurrent extension of $(X,\mathbb{P})$.

In [16], we provided necessary and sufficient conditions on the underlying Lévy process for the existence of recurrent extensions of $(X,\mathbb{P})$ that leave 0 by a jump. For the sake of completeness, we include an improved version of that result.

**Theorem 1.** *Let $(X,\mathbb{P})$ be an $\alpha$-self-similar Markov process that hits the cemetery point 0 in a finite time a.s. and $(\xi,\mathrm{P})$ the Lévy process associated with it via Lamperti's transformation. For $0 < \beta < 1/\alpha$, the following are equivalent:*

  (i) $\mathbf{E}(e^{\beta\xi_1}, 1 < \zeta) < 1$;
  (ii) $\mathbf{E}(I^{\alpha\beta}) < \infty$;
  (iii) *There exists a recurrent extension of $(X,\mathbb{P})$, say $X^{(\beta)}$, that leaves 0 by a jump and whose associated excursion measure $n^\beta$ is such that*

$$n^\beta(X_{0+} \in \mathrm{d}x) = c_{\alpha,\beta}\beta x^{-1-\beta}\,\mathrm{d}x, \qquad x > 0,$$

*where $c_{\alpha,\beta}$ is a constant.*

*In this case, the process $X^{(\beta)}$ is the unique recurrent extension of $(X,\mathbb{P})$ that leaves 0 by a jump distributed as $c_{\alpha,\beta}\eta_\beta$.*

The equivalence between (ii) and (iii) in Theorem 1 is the content of Proposition 1 in [16] and the equivalence between (i) and (ii) is a consequence of Lemma 2 below.

Thus, only the existence of recurrent extensions that leave 0 continuously remains to be established. In this vein, we proved in [16] that under the hypotheses:

(H2a) $(\xi,\mathbf{P})$ is non-arithmetic, that is, its state space is not a subgroup of $r\mathbb{Z}$, for any $r \in \mathbb{R}$,

(H2b) **Cramér's condition is satisfied**, that is, there exists a $\theta > 0$ s.t.

$$\mathbf{E}(e^{\theta\xi_1}, 1 < \zeta) = 1,$$

(H2c) for $\theta$ as in hypothesis (H2b), $\mathbf{E}(\xi_1^+ e^{\theta\xi_1}, 1 < \zeta) < \infty$

and provided $0 < \alpha\theta < 1$, there exists a recurrent extension of $(X,\mathbb{P})$ that leaves 0 continuously. In a previous work, Vuolle-Apiala [17] provided a sufficient condition on the resolvent of $(X,\mathbb{P})$ for the existence of recurrent extensions of $(X,\mathbb{P})$ that leave 0 continuously. Actually, in [16], we proved that in the case where the underlying Lévy process is non-arithmetic, the conditions of Vuolle-Apiala are equivalent to the conditions (H2b)–(H2c) above. So, it is natural to ask if the conditions of Vuolle-Apiala and those above are also necessary for the existence of recurrent extensions of $(X,\mathbb{P})$ that leave 0 continuously. The following counterexample answers this question negatively.

**Counterexample 1.** *Let $\sigma$ be a subordinator with law $P$ such that its law is not arithmetic and has some exponential moments of positive order, that is,*

$$\mathcal{E} := \{\lambda > 0,\ 1 < E(\mathrm{e}^{\lambda \sigma_1}) < \infty\} \neq \varnothing.$$

*Assume that the upper bound of $\mathcal{E}$, say $q$, belongs to $\mathcal{E} \cap ]0,1[$ and that the function*

$$m(x) := \mathbf{E}(1_{\{\sigma_1 > x\}} \mathrm{e}^{q\sigma_1}), \qquad x > 0,$$

*is regularly varying at infinity with index $-\beta$, for some $\beta \in ]1/2, 1[$. Let $(\xi, \mathbf{P})$ be the Lévy process with finite lifetime $\zeta$, obtained by killing $\sigma$ at an independent exponential time of parameter $\kappa = \log(E(\mathrm{e}^{q\sigma_1}))$. By construction, it follows that Cramér's condition*

$$\mathbf{E}(\mathrm{e}^{q\xi_1}, 1 < \zeta) = 1$$

*is satisfied and, by Karamata's theorem, the function*

$$m^\natural(x) := \int_0^x \mathbf{E}(1_{\{\xi_1 > u\}} \mathrm{e}^{q\xi_1}, 1 < \zeta)\, \mathrm{d}u, \qquad x \geq 0,$$

*is regularly varying at infinity with index $1 - \beta$. As a consequence, the integral $\mathbf{E}(\xi_1^+ \mathrm{e}^{q\xi_1}, 1 < \zeta)$ is not finite. We will denote by $\mathbf{P}^\natural$ the Girsanov-type transformation of $\mathbf{P}$ via the martingale $(\mathrm{e}^{q\xi_s}, s \geq 0)$, namely $\mathbf{P}^\natural$ is the unique measure s.t.*

$$\mathbf{P}^\natural = \mathrm{e}^{q\xi_t} \mathbf{P}, \qquad on\ \mathcal{F}_t,\ t \geq 0.$$

*Let $(X, \mathbb{P})$ be the 1-pssMp associated with $(\xi, \mathbf{P})$ via Lamperti's transformation and let $V_\lambda$ denote its $\lambda$-resolvent, $\lambda > 0$. We claim that the following assertions are satisfied:*

(P1) *for any $\lambda > 0$,*

$$\lim_{x \to 0+} m^\natural(\log(1/x)) \frac{V_\lambda f(x)}{x^q}$$
$$= \frac{1}{\Gamma(\beta)\Gamma(1-\beta)} \int_0^\infty f(y) \mathbf{E}^\natural\left(\exp\left\{-\lambda y \int_0^\infty \mathrm{e}^{-\xi_s}\,\mathrm{d}s\right\}\right) y^{-q}\,\mathrm{d}y$$

*for every $f: ]0, \infty[ \to \mathbb{R}$, continuous and with compact support;*

(P2) *the limit*

$$\lim_{x \to 0+} m^\natural(\log(1/x)) \frac{\mathbb{E}_x(1 - e^{-T_0})}{x^q} := C_q$$

*exists and $C_q \in ]0, \infty[$;*

(P3) *there exists a recurrent extension of $(X, \mathbb{P})$ that leaves $0$ continuously.*

That the properties (P1)–(P3) are satisfied in the framework of Counterexample 1 will be proven in Section 3.

In the previous counterexample, we have constructed a pssMp that satisfies neither the hypotheses of Vuolle-Apiala [17] nor all of the hypotheses in [16], but nonetheless admits a recurrent extension that leaves 0 continuously. This allowed us to realize that only Cramér's condition is relevant to the existence of recurrent extensions of pssMp. That is the content of the main theorem of this paper. To state the result, we need further notation.

First, observe that if Cramér's condition is satisfied with index $\theta$, then the process $M := (e^{\theta \xi_t}, t \geq 0)$ is a martingale under $\mathbf{P}$. In this case, we will denote by $\mathbf{P}^\natural$ the Girsanov-type transform of $\mathbf{P}$ via the martingale $M$, as we did in Counterexample 1. Under the law $\mathbf{P}^\natural$, the process $\xi$ is an $\mathbb{R}$-valued Lévy process with infinite lifetime and that drifts to $\infty$. We will denote by $J$ the exponential functional

$$J := \int_0^\infty \exp\{-\xi_s/\alpha\} \, ds.$$

A direct consequence of Theorem 1 of [4] and the fact that $(\xi, \mathbf{P}^\natural)$ is a Lévy process that drifts to $\infty$ is that $J < \infty$, $\mathbf{P}^\natural$-a.s. More details on the construction of the probability measure $\mathbf{P}^\natural$ and its properties can be found in Section 2.3 of [16]. We now have all of the elements necessary to state our main result.

**Theorem 2.** *Let $(X, \mathbb{P})$ be an $\alpha$-self-similar Markov process that hits its cemetery state $0$ in a finite time $\mathbb{P}$-a.s. and $(\xi, \mathbf{P})$ be the Lévy process associated with $(X, \mathbb{P})$ via Lamperti's transformation. The following are equivalent:*

(i) *there exists a $0 < \theta < 1/\alpha$ such that $\mathbf{E}(e^{\theta \xi_1}, 1 < \zeta) = 1$;*

(ii) *there exists a recurrent extension of $(X, \mathbb{P})$ that leaves $0$ continuously and such that its associated excursion measure from $0$, say $\mathbf{n}$, is such that*

$$\mathbf{n}(1 - e^{-T_0}) = 1.$$

*In this case, the recurrent extension in (ii) is unique and the entrance law associated with the excursion measure $\mathbf{n}$ is, for any $f$ positive and measurable, given by*

$$\mathbf{n}(f(X_t), t < T_0) = \frac{1}{t^{\alpha \theta} \Gamma(1 - \alpha \theta) \mathbf{E}^\natural(J^{\alpha \theta - 1})} \mathbf{E}^\natural\left(f\left(\frac{t^\alpha}{J^\alpha}\right) J^{\alpha \theta - 1}\right), \qquad t > 0, \qquad (3)$$

*with $\theta$ as in condition (i).*

Observe that condition (ii) of Theorem 2 implies that the inverse of the local time at 0 for the recurrent extension of $(X, \mathbb{P})$ is a stable subordinator of parameter $\alpha\theta$ for some $0 < \theta < 1/\alpha$. It is implicit in Theorem 2 that this is the unique $\theta > 0$ that fulfills condition (i), and vice versa. Moreover, the expression of the entrance law associated with **n** should be compared to the entrance law of Bertoin and Caballero [1] and Bertoin and Yor [3] for positive self-similar Markov processes that drift to $\infty$.

Besides, it is interesting to determine whether the recurrent extension in Theorem 2 is such that the underlying Lévy process satisfies the hypothesis (H2c) above. That is the content of the following corollary.

**Corollary 1.** *Assume that there exists a recurrent extension of $(X, \mathbb{P})$ that leaves $0$ continuously and let $\widetilde{\mathbb{P}}$ and **n** denote its law and excursion measure at $0$, respectively. For $\theta$ as in Theorem 2, the integrability condition*

$$\mathbf{E}(\xi_1^+ e^{\theta\xi_1}, 1 < \zeta) < \infty \tag{4}$$

*is satisfied if and only if*

$$\mathbf{n}(X_1^\theta, 1 < T_0) < \infty \tag{5}$$

*Furthermore, the latter holds if and only if*

$$\widetilde{\mathbb{E}}_x(X_t^\theta) < \infty \qquad \forall\, x \geq 0,\ \forall\, t \geq 0. \tag{6}$$

During the elaboration of this work, we learned that in [10], P. Fitzsimmons essentially proved the equivalence between (i) and (ii) in Theorem 2. He proved that Cramér's condition and a moment condition for the exponential functional $I$ are necessary and sufficient for the existence of a recurrent extension of $(X, \mathbb{P})$ that leaves $0$ continuously. Actually, the moment condition of Fitzsimmons is a consequence of Cramér's condition, as is proved in Lemma 2 below. Besides, Fitzsimmons' arguments and our own are completely different. He used arguments based on the theory of Kuznetsov measures and time-changes of processes with random birth and death, while our proof uses some general results on the excursions of pssMp obtained in our previous work [16].

The rest of this note is organized as follows. Section 2 is mainly devoted to the proof of Theorem 2 and in Section 3, we establish the facts claimed in Counterexample 1.

## 2. Proofs

To undertake our task, we need some notation. The Laplace exponent of $(\xi, \mathbf{P})$ is the function $\psi : \mathbb{R} \to \mathbb{R} \cup \{\infty\}$ defined by

$$\mathbf{E}(e^{\lambda\xi_1}, 1 < \zeta) := e^{\psi(\lambda)}, \qquad \lambda \in \mathbb{R}.$$

Hölder's inequality implies that $\psi$ is a strictly convex function on the set $\mathcal{E} := \{\lambda \in \mathbb{R} : \psi(\lambda) < \infty\}$. So, if Cramér's condition is satisfied, then the equation $\psi(\lambda) = 0$, $\lambda > 0$, has a unique root that we will denote hereafter by $\theta$. Observe that $[0,\theta] \subseteq \mathcal{E}$, that $\psi$ is derivable from the right at $0$ and from the left at $\theta$ and that

$$\mathbf{E}(\xi_1, 1 < \zeta) = \psi'_+(0) \in [-\infty, 0[, \qquad \mathbf{E}(\xi_1 e^{\theta \xi_1}, 1 < \zeta) = \psi'_-(\theta) \in ]0, \infty].$$

Our first purpose is to prove that (i) and (ii) in Theorem 1 are equivalent and that in Theorem 2, (i) implies (ii). To achieve this, we will need the following lemma.

**Lemma 2.** *Let $(\xi, \mathbf{P})$ be a Lévy process and assume that there is a $\beta \in ]0, 1/\alpha[$ such that*

$$\mathbf{E}(e^{\beta \xi_1}, 1 < \zeta) \leq 1.$$

*We then have that*

$$\mathbf{E}(I^{\alpha\beta - 1}) < \infty.$$

*Furthermore,*

$$\mathbf{E}(e^{\beta \xi_1}, 1 < \zeta) < 1, \quad \text{if and only if} \quad \mathbf{E}(I^{\alpha\beta}) < \infty.$$

**Proof.** For $t > 0$, let $Q_t$ denote the random variable

$$Q_t := \int_0^t \exp\{\xi_u/\alpha\} 1_{\{u < \zeta\}} \, du.$$

The main argument of the proof uses the fact that $\mathbf{E}(Q_t^{\alpha\beta}) < \infty$ for all $t > 0$. Indeed, the strict convexity of the mapping $\lambda \to \mathbf{E}(e^{\lambda \xi_1}, 1 < \zeta)$ implies that for any $p > 1$, $\mathbf{E}(e^{(\beta/p)\xi_t}, t < \zeta) = e^{t\psi(\beta/p)} < 1$, $t > 0$. Thus, for $p > 1$, we have

$$\mathbf{E}(Q_t^{\alpha\beta}) \leq t^{\alpha\beta} \mathbf{E}\left[\sup_{0 < u \leq t} \{e^{\beta \xi_u} 1_{\{u < \zeta\}}\}\right]$$

$$= t^{\alpha\beta} \mathbf{E}\left[\left(\sup_{0 < u \leq t} \{e^{(\beta/p)\xi_u} 1_{\{u < \zeta\}}\}\right)^p\right]$$

$$\leq t^{\alpha\beta} \mathbf{E}\left[\left(\sup_{0 < u \leq t} \{e^{(\beta/p)\xi_u} e^{-u\psi(\beta/p)} 1_{\{u < \zeta\}}\}\right)^p\right]$$

$$\leq t^{\alpha\beta} \left(\frac{p}{p-1}\right)^p \mathbf{E}\left[\{e^{(\beta/p)\xi_t} e^{-t\psi(\beta/p)} 1_{\{t < \zeta\}}\}^p\right]$$

$$\leq t^{\alpha\beta} \left(\frac{p}{p-1}\right)^p e^{-tp\psi(\beta/p)},$$

using Doob's $L_p$ inequality and the fact that the process $e^{(\beta/p)\xi_u - u\psi(\beta/p)}$, $u \geq 0$, is a positive martingale. We now prove the first claim in Lemma 2. On one hand, using the well-known inequality

$$||x|^{\alpha\beta} - |y|^{\alpha\beta}| \leq |x - y|^{\alpha\beta}, \qquad x, y \in \mathbb{R},$$

we get that

$$\mathbf{E}\left[\left(\int_0^\infty \exp\{\xi_s/\alpha\}1_{\{s<\zeta\}}\,\mathrm{d}s\right)^{\alpha\beta} - \left(\int_t^\infty \exp\{\xi_s/\alpha\}1_{\{s<\zeta\}}\,\mathrm{d}s\right)^{\alpha\beta}\right] \leq \mathbf{E}(Q_t^{\alpha\beta}) < \infty.$$

On the other hand, we have a.s.

$$\left(\int_0^\infty \exp\{\xi_s/\alpha\}1_{\{s<\zeta\}}\,\mathrm{d}s\right)^{\alpha\beta} - \left(\int_t^\infty \exp\{\xi_s/\alpha\}1_{\{s<\zeta\}}\,\mathrm{d}s\right)^{\alpha\beta}$$

$$= \alpha\beta \int_0^t \exp\{\xi_u/\alpha\}1_{\{u<\zeta\}}\left(\int_u^\infty \exp\{\xi_s/\alpha\}1_{\{s<\zeta\}}\,\mathrm{d}s\right)^{\alpha\beta-1}\mathrm{d}u$$

$$= \alpha\beta \int_0^t \exp\{\beta\xi_u\}1_{\{u<\zeta\}}\left(\int_0^\infty \exp\{\widetilde{\xi}_r/\alpha\}1_{\{r<\widetilde{\zeta}\}}\,\mathrm{d}r\right)^{\alpha\beta-1}\mathrm{d}u,$$

where $\widetilde{\xi}_r = \xi_{r+u} - \xi_u$, $r \geq 0$ and $\widetilde{\zeta} = \zeta - u$. Thus, by taking expectations, using Fubini's theorem and the independence of the increments of $\xi$, we obtain the identity

$$\mathbf{E}\left(\left(\int_0^\infty \exp\{\xi_s/\alpha\}1_{\{s<\zeta\}}\,\mathrm{d}s\right)^{\alpha\beta} - \left(\int_t^\infty \exp\{\xi_s/\alpha\}1_{\{s<\zeta\}}\,\mathrm{d}s\right)^{\alpha\beta}\right)$$

$$= \alpha\beta \int_0^t \mathbf{E}\left(\exp\{\beta\xi_u\}1_{\{u<\zeta\}}\left(\int_0^\infty \exp\{\widetilde{\xi}_r/\alpha\}1_{\{r<\widetilde{\zeta}\}}\,\mathrm{d}r\right)^{\alpha\beta-1}\right)\mathrm{d}u$$

$$= \alpha\beta\mathbf{E}(I^{\alpha\beta-1}) \int_0^t \mathbf{E}(\exp\{\beta\xi_u\}1_{\{u<\zeta\}})\,\mathrm{d}u.$$

The first claim in Lemma 2 follows. To prove the second assertion, we first assume that $\mathbf{E}(\mathrm{e}^{\beta\xi_1}, 1 < \zeta) < 1$. Thus, by letting $t$ tend to infinity and integrating in the latter equation, we obtain the identity

$$\mathbf{E}(I^{\alpha\beta}) = \frac{\alpha\beta}{\psi(\beta)}\mathbf{E}(I^{\alpha\beta-1}). \tag{7}$$

This relation is well known; see, for example, [4] and [15]. Together with the first assertion of the lemma, this implies that $\mathbf{E}(I^{\alpha\beta}) < \infty$. We now prove the reciprocal. If $\mathbf{E}(I^{\alpha\beta}) < \infty$, then we have

$$\infty > \mathbf{E}\left(\left(\int_0^\zeta \exp\{\xi_s/\alpha\}\,\mathrm{d}s\right)^{\alpha\beta}\right)$$

$$> \mathbf{E}\left(\left(\int_1^\zeta \exp\{\xi_s/\alpha\}\,\mathrm{d}s\right)^{\alpha\beta}1_{\{1<\zeta\}}\right) \tag{8}$$

$$= \mathbf{E}\left(\mathrm{e}^{\beta\xi_1}\mathbf{E}\left(\left(\int_0^\infty \exp\{(\xi_{1+s}-\xi_1)/\alpha\}1_{\{1+s<\zeta\}}\,\mathrm{d}s\right)^{\alpha\beta}\right)1_{\{1<\zeta\}}\right)$$

$$= \mathbf{E}(\mathrm{e}^{\beta\xi_1}1_{\{1<\zeta\}})\mathbf{E}\bigg(\bigg(\int_0^\zeta \exp\{\xi_s/\alpha\}\,\mathrm{d}s\bigg)^{\alpha\beta}\bigg),$$

due to the fact that $\xi$ is a Lévy process. So, we have that, in this case, $\mathbf{E}(\mathrm{e}^{\beta\xi_1}, 1<\zeta)<1$. □

**Proof of Theorem 2: (i) implies (ii).** The proof of this is based on Theorem 3 of [16], but to use that result, we first need to establish some weak duality relations.

By assumption (i) and Lemma 2, we have that $\mathbf{E}(I^{\alpha\theta-1})<\infty$. Moreover, let $(\xi, \widehat{\mathbf{P}}^\natural) := (-\xi, \mathbf{P}^\natural)$ denote the dual of $(\xi, \mathbf{P}^\natural)$. Then, $(\xi, \widehat{\mathbf{P}}^\natural)$ drifts to $-\infty$ because $(\xi, \mathbf{P}^\natural)$ drifts to $\infty$ and, as a consequence, $I<\infty$, $\widehat{\mathbf{P}}^\natural$-a.s. Furthermore, $(\xi, \widehat{\mathbf{P}}^\natural)$ satisfies the hypotheses of Lemma 2 with $\beta = \theta$, due to the identity

$$\widehat{\mathbf{E}}^\natural(\mathrm{e}^{\theta\xi_1}) = \mathbf{E}^\natural(\mathrm{e}^{-\theta\xi_1}) = \mathbf{E}(\mathrm{e}^{-\theta\xi_1}\mathrm{e}^{\theta\xi_1}, 1<\zeta)\le 1.$$

Thus, we can also ensure that $\widehat{\mathbf{E}}^\natural(I^{\alpha\theta-1})<\infty$. Now, let $\widehat{\mathbb{P}}^\natural$ be the law of the $\alpha$-pssMp associated with $(\xi, \widehat{\mathbf{P}}^\natural)$ via Lamperti's transformation. $(X, \widehat{\mathbb{P}}^\natural)$ is then an $\alpha$-pssMp that hits 0 continuously and in a finite time, $\widehat{\mathbb{P}}^\natural$-a.s., and, according to Lemma 2 in [3], $(X, \mathbb{P}^\natural)$ and $(X, \widehat{\mathbb{P}}^\natural)$ are in weak duality with respect to the measure $\alpha^{-1}x^{1/\alpha-1}\,\mathrm{d}x, x>0$, and, given that the law $\mathbb{P}^\natural$ is the $h$-transform of the law $\mathbb{P}$ via the invariant function $h(x)=x^\theta$ for the semigroup of $(X, \mathbb{P})$ (see Proposition 5 of [16]), it then follows that $(X, \mathbb{P})$ and $(X, \widehat{\mathbb{P}}^\natural)$ are in weak duality w.r.t. the measure $\alpha^{-1}x^{1/\alpha-1-\theta}\,\mathrm{d}x$, $x>0$. Furthermore, we have that for any $\lambda>0$,

$$\alpha^{-1}\int_0^\infty \mathrm{d}x\, x^{1/\alpha-1-\theta}\mathbb{E}_x(\mathrm{e}^{-\lambda T_0})<\infty, \qquad \alpha^{-1}\int_0^\infty \mathrm{d}x\, x^{1/\alpha-1-\theta}\widehat{\mathbb{E}}^\natural_x(\mathrm{e}^{-\lambda T_0})<\infty. \quad (9)$$

Indeed, for $\lambda>0$,

$$\alpha^{-1}\int_0^\infty \mathrm{d}x\, x^{1/\alpha-1-\theta}\mathbb{E}_x(\mathrm{e}^{-\lambda T_0}) = \alpha^{-1}\int_0^\infty \mathrm{d}x\, x^{1/\alpha-1-\theta}\mathbf{E}(\mathrm{e}^{-\lambda x^{1/\alpha}I})$$

$$= \mathbf{E}\bigg(\alpha^{-1}\int_0^\infty \mathrm{d}x\, x^{1/\alpha-1-\theta}\mathrm{e}^{-\lambda x^{1/\alpha}I}\bigg)$$

$$= \lambda^{\alpha\theta-1}\mathbf{E}(I^{\alpha\theta-1})\Gamma(1-\alpha\theta)<\infty.$$

The same calculation applies to the verification of the finiteness of the second integral in equation (9). This being said, Theorem 3 of [16] ensures that there exists a unique recurrent extension of $(X, \mathbb{P})$ such that the $\lambda$-resolvent of its excursion measure, say $\mathbf{n}$, is given by

$$\mathbf{n}\bigg(\int_0^{T_0}\mathrm{e}^{-\lambda t}f(X_t)\,\mathrm{d}t\bigg) = \frac{1}{\alpha\Gamma(1-\alpha\theta)\widehat{\mathbf{E}}^\natural(I^{\alpha\theta-1})}\int_0^\infty f(x)x^{1/\alpha-1-\theta}\widehat{\mathbb{E}}^\natural_x(\mathrm{e}^{-\lambda T_0})\,\mathrm{d}x, \quad (10)$$

for $\lambda\ge 0$, and any function $f$, positive and measurable on $[0,\infty[$. An easy calculation proves that the $\lambda$-resolvent of $\mathbf{n}$ satisfies the self-similarity property in Lemma 2 of [16] and therefore the excursion measure **n** is self-similar. In particular, $\mathbf{n}(1 - e^{-T_0}) = 1$ and the potential of **n** is given by

$$\mathbf{n}\left(\int_0^{T_0} f(X_t)\,dt\right) = \frac{1}{\alpha\Gamma(1-\alpha\theta)\widehat{\mathbf{E}}^\natural(I^{\alpha\theta-1})} \int_0^\infty f(x) x^{1/\alpha-1-\theta}\,dx.$$

Compared with the result in Lemma 3 of [16], this implies that

$$\widehat{\mathbf{E}}^\natural(I^{\alpha\theta-1}) = \mathbf{E}(I^{\alpha\theta-1}). \tag{11}$$

Actually, Theorem 3 of [16] also establishes that there exists a recurrent extension of $(X, \widehat{\mathbb{P}}^\natural)$ with excursion measure $\widehat{\mathbf{n}}$ such that

$$\widehat{\mathbf{n}}\left(\int_0^{T_0} e^{-\lambda t} f(X_t)\,dt\right) = \frac{1}{\alpha\Gamma(1-\alpha\theta)\mathbf{E}(I^{\alpha\theta-1})} \int_0^\infty f(x) x^{1/\alpha-1-\theta} \mathbb{E}_x(e^{-\lambda T_0})\,dx.$$

Moreover, the recurrent extensions of $(X, \mathbb{P})$ and $(X, \widehat{\mathbb{P}}^\natural)$ associated with **n** and $\widehat{\mathbf{n}}$, respectively, are still in weak duality. To verify that **n** is carried by the paths that leave 0 continuously, we claim that the image of **n** under time reversal at time $T_0$ is $\widehat{\mathbf{n}}$. This follows from the fact that **n** and $\widehat{\mathbf{n}}$ have the same potential and an application of a result for time reversal of Kuznetsov measures established in Section XIX.23 of Dellacherie, Maisonneuve and Meyer [8]. Thus, using the Markov property and the fact that $(X, \widehat{\mathbb{P}}^\natural)$ is a pssMp that hits 0 continuously and in a finite time $\widehat{\mathbb{P}}^\natural$-a.s., given that the underlying Lévy process $(\xi, \widehat{\mathbf{P}}^\natural)$ drifts to $-\infty$, we get that $\widehat{\mathbf{n}}$ is carried by the paths that hit 0 continuously and therefore

$$0 = \widehat{\mathbf{n}}(X_{T_0-} > 0) = \mathbf{n}(X_{0+} > 0).$$

□

**Proof of Theorem 2: (ii) implies (i).** Assume that the hypothesis (ii) in Theorem 2 holds and denote by $\widetilde{X}$ a recurrent extension of $(X, \mathbb{P})$ with excursion measure **n**. We claim that, in this case, we have the inequality

$$\mathbf{E}(e^{\vartheta\xi_1}, 1 < \zeta) \leq 1 \tag{12}$$

and that, in fact, the strict inequality is impossible, that is, Cramér's condition is satisfied. Taking for granted inequality (12), it is easy to see that the latter holds. Assume that $\mathbf{E}(e^{\vartheta\xi_1}, 1 < \zeta) < 1$. Theorem 1 would imply that $(X, \mathbb{P})$ admits a recurrent extension that leaves 0 by a jump and with jumping-in measure proportional to $\eta_\vartheta$. This implies that the measure $m = 2^{-1}\mathbf{n} + 2^{-1}c_{\alpha,\vartheta}\mathbb{P}_{\eta_\vartheta}$ is a self-similar excursion measure compatible with $(X, \mathbb{P})$ and with index of self-similarity $\alpha\vartheta$; as before, $c_{\alpha,\vartheta}$ is a normalizing constant. Therefore, there exists a recurrent extension of $(X, \mathbb{P})$ with excursion measure $m$ that may leave 0 by a jump and continuously, which leads to a contradiction of the fact that any recurrent extension of $(X, \mathbb{P})$ either leaves 0 by a jump or continuously.

We will now prove that inequality (12) holds. It follows from Lemmas 2 and 3 of [16] that there exists a $\vartheta \in {]}0, 1/\alpha[$ such that the potential of the measure $\mathbf{n}$ is given by

$$\nu(\mathrm{d}y) := \mathbf{n}\left(\int_0^{T_0} 1_{\{X_t \in \mathrm{d}y\}}\,\mathrm{d}t\right) = C_{\alpha,\alpha\vartheta} y^{(1-\alpha-\alpha\vartheta)/\alpha}\,\mathrm{d}y, \qquad y > 0,$$

for a constant $0 < C_{\alpha,\vartheta} < \infty$ and that $\nu$ is the unique invariant measure for $\widetilde{X}$; the uniqueness holds up to a multiplicative constant. It follows that $\nu$ is an excessive measure for $(X, \mathbb{P})$. Besides, the Revuz measure of the additive functional $B$ defined by $B_t := \int_0^t X_s^{-1/\alpha}\,\mathrm{d}s$, $0 \leq t < T_0$, relative to $\nu$, is given by

$$\nu_B(\mathrm{d}y) := C_{\alpha,\alpha\vartheta} y^{-1-\vartheta}\,\mathrm{d}y, \qquad y > 0.$$

So, due to Lamperti's transformation, the process $(\mathrm{e}^\xi, \mathrm{P})$ is obtained by time-changing $(X, \mathbb{P})$ by the right-continuous inverse of $B$, thus the measure $\nu_B$ is excessive for $(\mathrm{e}^\xi, \mathbf{P})$. It follows from this that the measure $\mathrm{e}^{-\vartheta y}\,\mathrm{d}y$, $y \in \mathbb{R}$, is excessive for $(\xi, \mathbf{P})$. The latter assertion implies that for every positive and bounded function $f : \mathbb{R} \to \mathbb{R}$, the following inequalities hold:

$$\int_{\mathbb{R}} \mathrm{e}^{-\vartheta x} f(x)\,\mathrm{d}x \geq \int_{\mathbb{R}} \mathrm{e}^{-\vartheta x} \mathbf{E}_x(f(\xi_1), 1 < \zeta)\,\mathrm{d}x$$

$$= \mathbf{E}\left(\int_{\mathbb{R}} \mathrm{e}^{-\vartheta(y-\xi_1)} f(y)\,\mathrm{d}y, 1 < \zeta\right)$$

$$= \int_{\mathbb{R}} \mathbf{E}(\mathrm{e}^{\vartheta \xi_1}, 1 < \zeta) f(y) \mathrm{e}^{-\theta y}\,\mathrm{d}y.$$

From these, inequality (12) follows. $\square$

We have thus completed the proof of the equivalence between assertions (i) and (ii) in Theorem 2. Observe that the $\theta$ in the proof of the implication (i) $\Longrightarrow$ (ii) is equal to the $\vartheta$ in the implication (ii) $\Longrightarrow$ (i).

We next prove the uniqueness and characterization of the entrance law associated with the excursion measure claimed in Theorem 2.

### 2.1. Uniqueness and characterization

Assume that there exist two recurrent extensions of $(X, \mathbb{P})$ that satisfy the conditions of (ii) of Theorem 2 and let $\mathbf{n}$ and $\mathbf{n}'$ be its associated excursions measures. There then exist $\theta_1$ and $\theta_2$ such that Cramér's condition is satisfied. The strict convexity of the mapping $\lambda \to \mathbb{E}(\mathrm{e}^{\lambda \xi_1}, 1 < \zeta)$ implies that $\theta_2 = \theta = \theta_1$. As a consequence, the potential of both excursion measures is given by equation (2), with $\gamma$ replaced by $\alpha\theta$. Therefore, arguing as in the proof of Lemma 1, we show that $\mathbf{n} = \mathbf{n}'$, which completes the proof of uniqueness.

The characterization of the entrance law follows from our proof of the fact that (i) implies (ii) in Theorem 2. On one hand, by construction, the resolvent of the excursion measure **n** is given by equation (10). On the other hand,

$$\mathbf{n}\left(\int_0^{T_0} e^{-\lambda t} f(X_t)\,dt\right) = \int_0^\infty e^{-\lambda t} t^{-\alpha\theta} \mathbf{n}(f(t^\alpha X_1), 1 < T_0)\,dt$$
$$= \mathbf{n}\left(\int_0^\infty du\,\alpha^{-1} u^{1/\alpha-1-\theta} f(u) X_1^{\theta-1/\alpha} e^{-\lambda u^{1/\alpha} X_1^{-1/\alpha}} 1_{\{1<T_0\}}\right) \quad (13)$$
$$= \int_0^\infty du\,\alpha^{-1} u^{1/\alpha-1-\theta} f(u) \mathbf{n}(X_1^{\theta-1/\alpha} e^{-\lambda u^{1/\alpha} X_1^{-1/\alpha}} 1_{\{1<T_0\}}),$$

where we used Fubini's theorem three times, combined with the scaling property of **n** and a change of variables. Comparing the results in equations (10) and (13), we get the identity

$$\mathbf{n}(X_1^{\theta-1/\alpha} \exp\{-\lambda u^{1/\alpha} X_1^{-1/\alpha}\} 1_{\{1<T_0\}}) = \frac{1}{\Gamma(1-\alpha\theta)\widehat{\mathbf{E}}^\natural(I^{\alpha\theta-1})} \widehat{\mathbb{E}}_u^\natural(e^{-\lambda T_0})$$

for all $\lambda \geq 0$ and a.e. $u > 0$. As a consequence,

$$\mathbf{n}(X_1^{\theta-1/\alpha} 1_{\{1<T_0\}}) < \infty.$$

By the dominated convergence theorem, the latter identity holds for all $\lambda \geq 0$ and all $u > 0$. Recall that, by Lamperti's transformation, $T_0$ under $\widehat{\mathbb{P}}_u^\natural$ has the same law as $u^{1/\alpha} I$ under $\widehat{\mathbf{P}}^\natural$. So, by the uniqueness of Laplace transforms, it follows that

$$\mathbf{n}(X_1^{\theta-1/\alpha} f(X_1^{-1/\alpha}) 1_{\{1<T_0\}}) = \frac{1}{\Gamma(1-\alpha\theta)\widehat{\mathbf{E}}^\natural(I^{\alpha\theta-1})} \widehat{\mathbf{E}}^\natural(f(I)).$$

The claim in Theorem 2 follows from this identity using the scaling property of **n** and the fact that the law of $I$ under $\widehat{\mathbf{P}}^\natural$ is equal to that of $J$ under $\mathbf{P}^\natural$.

Having completed the proof of Theorem 2, we next prove Corollary 1.

## 2.2. Proof of Corollary 1

Due to the existence of the left derivative of $\psi$ at $\theta$, it follows from Proposition 3.1 of [7] that

$$\mathbf{E}^\natural(J^{-1}) = \widehat{\mathbf{E}}^\natural(I^{-1}) = -\widehat{\mathbf{E}}^\natural(\xi_1) = \mathbf{E}(\xi_1 e^{\theta\xi_1}, 1 < \zeta)$$

and the leftmost quantity is finite if and only if $\mathbf{E}(\xi_1^+ e^{\theta\xi_1}, 1 < \zeta) < \infty$. So, that (4) and (5) are equivalent is an easy consequence of the representation of the entrance law obtained in Theorem 2. We will next prove that (5) is equivalent to (6). Let $V_\lambda$ denote the $\lambda$-resolvent of $(X, \mathbb{P})$ and $U_\lambda$ be the $\lambda$-resolvent of the unique recurrent extension of $(X, \mathbb{P})$ that leaves 0 continuously. The invariance of $h(x) = x^\theta, x > 0$, implies that

$$\mathbf{n}(X_1^\theta, 1 < T_0) = \mathbf{n}(X_t^\theta, t < T_0), \qquad t > 0.$$

Using a well-known decomposition formula and Fubini's theorem, we get

$$\begin{aligned}
\lambda U_\lambda h(x) &= \lambda V_\lambda h(x) + \mathbb{E}_x(e^{-\lambda T_0}) \lambda U_\lambda h(0) \\
&= h(x) + \mathbb{E}_x(e^{-\lambda T_0}) \frac{\lambda \mathbf{n}(\int_0^{T_0} e^{-\lambda t} h(X_t) \, dt)}{\mathbf{n}(1 - e^{-\lambda T_0})} \\
&= h(x) + \mathbb{E}_x(e^{-\lambda T_0}) \lambda^{-\alpha \theta} \int_0^\infty \lambda e^{-\lambda t} \mathbf{n}(h(X_t), t < T_0) \, dt \\
&= h(x) + \mathbb{E}_x(e^{-\lambda T_0}) \lambda^{-\alpha \theta} \mathbf{n}(h(X_1), 1 < T_0).
\end{aligned}$$

Thus, $\lambda U_\lambda h(x) < \infty$ for all $x$ if and only if $\mathbf{n}(X_1^\theta, 1 < T_0) < \infty$. From this, we get that if (5) holds, then $\widetilde{\mathbb{E}}_x(X_t^\theta) < \infty$ for all $x > 0$ and a.e. $t > 0$. The self-similarity implies that, in this case, the latter holds for all $x > 0$ and all $t > 0$. This completes the proof of Corollary 1.

## 3. Proof of Counterexample 1

A key tool in the establishment of (P1) and (P2) is the following version of Erickson's renewal theorem [9].

**Lemma 3** (Erickson's renewal theorem [9]). *Let $G$ be a non-arithmetic probability distribution function on $\mathbb{R}^+$ such that $1 - G$ is a regularly varying function at infinity with index $\gamma \in \,]1/2, 1]$, $U$ the renewal measure associated with $G$ and $m(x) := \int_0^x (1 - G(u)) \, du$, $x \geq 0$. Then:*

(i) *for any directly Riemann integrable function $g: \mathbb{R}^+ \to \mathbb{R}^+$,*

$$\lim_{t \to \infty} m(t) \int_0^t g(t - y) U(dy) = \frac{1}{\Gamma(\gamma)\Gamma(1 - \gamma)} \int_0^\infty g(y) \, dy;$$

(ii) *for any directly Riemann integrable function $g: \mathbb{R} \to \mathbb{R}^+$,*

$$\lim_{t \to \infty} m(t) \int_{-\infty}^\infty g(y - t) U(dy) = \frac{1}{\Gamma(\gamma)\Gamma(1 - \gamma)} \int_{-\infty}^\infty g(y) \, dy.$$

The statement in (i) of Lemma 3 is the content of Erickson's renewal Theorem 3 and so only (ii) requires a proof, which is postponed to the end of this section. Next, we proceed to prove the claims in Counterexample 1. To that end, observe that the law $\mathbf{P}^\natural$ is that of a subordinator with infinite lifetime, such that the tail probability $\mathbf{P}^\natural(\xi_1 > x)$ is a regularly varying function with index $\beta \in \,]1/2, 1[$. Let $U^\natural$ be the renewal measure of the subordinator with law $\mathbf{P}^\natural$, that is,

$$U^\natural(\mathrm{d}y) = \int_0^\infty \mathbf{P}^\natural(\xi_t \in \mathrm{d}y)\,\mathrm{d}t, \qquad y \geq 0.$$

According to Bertoin and Doney [2], the measure $U^\natural$ is the renewal measure associated with the probability distribution function given by $F(\cdot) = \mathbf{P}^\natural(\xi_\mathrm{e} \leq \cdot)$, where e is a standard exponential r.v. independent of $\xi$ under $\mathbf{P}^\natural$. Let $\mathbb{P}^\natural$ be the law of the 1-pssMp associated with $(\xi, \mathbf{P}^\natural)$ via Lamperti's transformation. The measure $\mathbb{P}^\natural$ is such that

$$\mathbb{P}^\natural = X_t^q \mathbb{P} \qquad \text{on } \mathcal{G}_t,\ t \geq 0.$$

It follows that the resolvents of $(X, \mathbb{P}^\natural)$ and $(X, \mathbb{P})$ are related by

$$V_\lambda^\natural f(x) = \frac{V_\lambda f h_q(x)}{h_q(x)}, \qquad x \in\,]0, \infty[, \tag{14}$$

with $h_q(x) := x^q, x > 0$. Moreover, we have that, for any function $f: \mathbb{R}^+ \to \mathbb{R}^+$ such that the mapping $y \to f(\mathrm{e}^y)\mathrm{e}^y$ is directly Riemann integrable,

$$\lim_{x \to 0+} m^\natural(\log(1/x)) V_0^\natural f(x) = \frac{1}{\Gamma(1-\gamma)\Gamma(\gamma)} \int_0^\infty f(y)\,\mathrm{d}y. \tag{15}$$

Indeed, by applying Lamperti's representation and (ii) in Lemma 3, we obtain

$$m^\natural(\log(1/x)) V_0^\natural f(x) = m^\natural(\log(1/x)) \mathbf{E}^\natural\left[\int_0^\infty f(x\mathrm{e}^{\xi_t}) x \mathrm{e}^{\xi_t}\,\mathrm{d}t\right]$$

$$= m^\natural(\log(1/x)) \int_\mathbb{R} f(\mathrm{e}^{y-\log(1/x)}) \mathrm{e}^{y-\log(1/x)} U^\natural(\mathrm{d}y)$$

$$\xrightarrow[x \to 0+]{} \frac{1}{\Gamma(1-\gamma)\Gamma(\gamma)} \int_\mathbb{R} f(\mathrm{e}^y) \mathrm{e}^y\,\mathrm{d}y$$

and by making a change of variables in the rightmost quantity, we obtain (15). Moreover, repeating the arguments at the beginning of Section 3 in [3], we prove that for every $f: ]0, \infty[ \to \mathbb{R}$, continuous and with compact support, and $\lambda > 0$,

$$\lim_{x \to 0+} m^\natural(\log(1/x)) V_\lambda^\natural f(x) = \frac{1}{\Gamma(1-\gamma)\Gamma(\gamma)} \int_0^\infty f(y) \mathbf{E}^\natural\left[\exp\left\{-\lambda y \int_0^\infty \mathrm{e}^{-\xi_s}\,\mathrm{d}s\right\}\right]\mathrm{d}y. \tag{16}$$

Therefore, the claim in (P1) is a straightforward consequence of (14) and (16).

Besides, in Lemma 4 of [16] we proved that, in general, the exponential functional $I$ satisfies the equation in law

$$I \stackrel{\mathrm{Law}}{=} Q + M\widetilde{I}, \qquad \text{with } (Q, M) := \left(\int_0^1 \exp\{\xi_s/\alpha\} 1_{\{s<\zeta\}}\,\mathrm{d}s, \mathrm{e}^{\alpha^{-1}\xi_1} 1_{\{1<\zeta\}}\right) \text{ and } I \stackrel{\mathrm{Law}}{=} \widetilde{I},$$

and the pair $(Q, M)$ is independent of $\widetilde{I}$. Moreover, under the hypotheses (H2) of Lemma 4 in [16], we obtained, as a consequence of Goldie's Theorems 2.3 and 4.1 in [13], an estimate of the tail probability of $I$. A perusal of the proofs provided by Goldie for those theorems allows us to ensure that the arguments can be extended, using Erickson's renewal theorem instead of the classical renewal theorem, to prove the following lemma.

**Lemma 4.** *Under the hypothesis of Counterexample 1, we have that*

$$\lim_{t \to \infty} m^\natural(\log(t)) t^q \mathbb{P}_1(T_0 > t)$$
$$= \frac{1}{\Gamma(1-\gamma)\Gamma(\gamma)} \mathbf{E}\left(\left(\int_0^\infty \exp\{\xi_s\} 1_{\{s<\zeta\}} \, ds\right)^q - \left(\int_1^\infty \exp\{\xi_s\} 1_{\{s<\zeta\}} \, ds\right)^q\right)$$
$$= \frac{1}{\Gamma(1-\gamma)\Gamma(\gamma)} q \mathbf{E}(I^{q-1}) \in \,]0, \infty[.$$

Therefore, Lemma 4 and Karamata's Tauberian theorem together imply that the property (P2) is satisfied.

*Remark 1.* The expression of the value of the limit in Lemma 4 is a consequence of the proof of Lemma 2.

Finally, to prove that condition (P3) is satisfied, we argue, as in [17], pages 556–557, to ensure that there exists a family of finite measures on $]0, \infty[$, say $(\mathbf{n}_\lambda, \lambda > 0)$, such that

$$\mathbf{n}_\lambda f = \lim_{x \to 0+} \frac{V_\lambda f(x)}{\mathbb{E}_x(1 - e^{-T_0})}$$

for any $f$, continuous and with compact support on $]0, \infty[$, and for $\lambda > 0$. Moreover, the family $(\mathbf{n}_\lambda, \lambda > 0)$ satisfies the resolvent-type equation for $\lambda, \mu > 0$,

$$\mathbf{n}_\lambda V_\mu f = \frac{\mathbf{n}_\mu f - \mathbf{n}_\lambda f}{\lambda - \mu},$$

for any $f$ continuous and with bounded support on $]0, \infty[$. Thus, Theorem 6.9 of [12] and Theorem 4.7 of [6] imply that there exist a unique excursion measure $\mathbf{n}$ whose $\lambda$-potential is equal to $\mathbf{n}_\lambda$,

$$\mathbf{n}\left(\int_0^{T_0} e^{-\lambda t} 1_{\{X_t \in dy\}} \, dt\right) = \mathbf{n}_\lambda(dy),$$

for any $\lambda > 0$. In fact, all of the results of Vuolle-Apiala [17] are still valid if we replace the power function that gives the normalization in his hypotheses (Aa) and (Ab) by a regularly varying function. Therefore, Theorem 1.2 of [17] ensures that $\mathbf{n}(X_{0+} > 0) = 0$. According to Blumenthal's theorem [5], associated with this excursion measure $\mathbf{n}$, there exists a unique recurrent extension of $(X, \mathbb{P})$ that leaves 0 continuously. This completes the proof of Counterexample 1. Now, we just have to prove that (ii) in Lemma 3 holds.

**Proof of Lemma 3.** The claim in (i) is Theorem 3 of Erickson [9] and that (ii) holds is a consequence of the latter. We next prove the result for step functions and the general case follows by a standard argument. Let $(a_k, k \in \mathbb{Z})$ be a sequence of positive real numbers such that $\sum_{k \in \mathbb{Z}} a_k < \infty$ and $h > 0$ a constant. A consequence of Theorem 1 of Erickson [9] is that for any $k \in \mathbb{N}$,

$$m(t+kh) \int_{\mathbb{R}} 1_{\{[kh,(k+1)h[\}}(y-t) U(\mathrm{d}y) \xrightarrow[t \to \infty]{} C_\gamma \int_{\mathbb{R}} 1_{\{[0,h[\}}(y) \, \mathrm{d}y = C_\gamma \int_{\mathbb{R}} 1_{\{[kh,(k+1)h[\}}(y) \, \mathrm{d}y,$$

with $C_\gamma = (\Gamma(\gamma)\Gamma(1-\gamma))^{-1}$, and uniformly in $k$. Thus, given that $m$ is an increasing function, we get that

$$m(t) \int_{\mathbb{R}} \sum_{k \in \mathbb{N}} a_k 1_{\{[kh,(k+1)h[\}}(y-t) U(\mathrm{d}y)$$

$$\leq \sum_{k \in \mathbb{N}} a_k m(t+kh) \int_{\mathbb{R}} 1_{\{[kh,(k+1)h[\}}(y-t) U(\mathrm{d}y).$$

Therefore,

$$\limsup_{t \to \infty} m(t) \int_{\mathbb{R}} \sum_{k \in \mathbb{N}} a_k 1_{\{[kh,(k+1)h[\}}(y-t) U(\mathrm{d}y) \leq C_\gamma \sum_{k \in \mathbb{N}} a_k \int_{\mathbb{R}} 1_{\{[kh,(k+1)h[\}}(y) \, \mathrm{d}y$$

$$\leq C_\gamma \int_{\mathbb{R}} \sum_{k \in \mathbb{N}} a_k 1_{\{[kh,(k+1)h[\}}(y) \, \mathrm{d}y.$$

Because $m$ is regularly varying with positive index, the following limit

$$\lim_{t \to \infty} \frac{m(t)}{m(t+kh)} = 1,$$

holds uniformly in $k \in \mathbb{N}$. A standard application of Fatou's theorem and an easy manipulation gives that

$$\liminf_{t \to \infty} m(t) \int_{\mathbb{R}} \sum_{k \in \mathbb{N}} a_k 1_{\{[kh,(k+1)h[\}}(y-t) U(\mathrm{d}y) \geq C_\gamma \int_{\mathbb{R}} \sum_{k \in \mathbb{N}} a_k 1_{\{[kh,(k+1)h[\}}(y) \, \mathrm{d}y.$$

Let $g$ be the step function defined by

$$g(t) = \sum_{k \in \mathbb{Z}} a_k 1_{[kh,(k+1)h[}(t), \qquad t \in \mathbb{R}.$$

It follows from the arguments above that

$$\lim_{t \to \infty} m(t) \int_{\mathbb{R}} g(y-t) 1_{\{y-t \geq 0\}} U(\mathrm{d}y) = C_\gamma \int_{\mathbb{R}} g(y) 1_{\{y \geq 0\}} \, \mathrm{d}y.$$

Moreover, the assertion in (i) implies that

$$\lim_{t\to\infty} m(t) \int_{\mathbb{R}} g(y-t) 1_{\{y-t<0\}} U(\mathrm{d}y) = \lim_{t\to\infty} m(t) \int_0^t g(-(t-y)) U(\mathrm{d}y) = C_\gamma \int_0^\infty g(-y)\,\mathrm{d}y,$$

from which the result follows. $\square$

## Acknowledgements

I would like to thank Professor P. Fitzsimmons for pointing out an error in an early draft of this article. I am indebted to the anonymous referee for insightful comments and for suggesting the present proof of inequality (12). Research was supported by a Research Award from the Council of Science and Technology of the state of Guanajuato, Mexico (CONCYTEG).